# Gauss Maps of the Mean Curvature Flow

Mu-Tao Wang

## Abstract

Let $F: \Sigma^n \times [0, T) \to \mathbb{R}^{n+m}$ be a family of compact immersed submanifolds moving by their mean curvature vectors. We show the Gauss maps $\gamma : (\Sigma^n, g_t) \to G(n, m)$ form a harmonic heat flow with respect to the time-dependent induced metric $g_t$. This provides a more systematic approach to investigating higher codimension mean curvature flows. A direct consequence is any convex function on $G(n, m)$ produces a subsolution of the nonlinear heat equation on $(\Sigma, g_t)$. We also show the condition that the image of the Gauss map lies in a totally geodesic submanifold of $G(n, m)$ is preserved by the mean curvature flow. Since the space of Lagrangian subspaces is totally geodesic in $G(n, n)$, this gives an alternative proof that any Lagrangian submanifold remains Lagrangian along the mean curvature flow.

## 1 Introduction

The maximum principle has been exploited by S-T. Yau and his coauthors to obtain gradient estimates for various geometric nonlinear partial differential equations since the '70. The principal idea is to identify a suitable expression of the gradient of the solution to which the differential operator (usually Laplace or heat operator) is applied. The prototype is the so-called Bochner formula and the geometry (curvature) of the manifold arises naturally in the calculation. Such techniques were extended by R. Hamilton to Ricci flows which are nonlinear parabolic systems. The curvature as second derivatives of the metric satisfies a certain parabolic system and the geometry of the space of curvature is used to derived estimates. A powerful theorem of Hamilton says any convex invariant subset of the space of curvature is preserved by the Ricci flow.

Such curvature estimates were extended to the mean curvature flow of hypersurfaces by G. Huisken. Recall a mean curvature flow is an evolution equation under which a submanifold deforms in the direction of its mean curvature vector. The mean curvature flow of an immersion $F : \Sigma^n \to \mathbb{R}^{n+m}$ is a family of immersions parametrized by $t$, $F : \Sigma \times [0, T) \to \mathbb{R}^{n+m}$ that satisfies

$$\frac{d}{dt}F(x,t) = H(x,t)$$
$$F_0 = F.$$

We shall denote the image $F(\cdot, t)$ by $\Sigma_t$, then $H(x,t)$ is the mean curvature vector of $\Sigma_t$ at $F(x,t)$.

Codimension-one mean curvature flows have been studied by Huisken et al and many beautiful results are obtained. In this paper, we shall focus on higher codimension mean curvature flows, i.e. the $m > 1$ case.

Consider a system with $n$ variables and $n + m$ functions. The configuration space of the differentials of the solutions consists of $n$ vectors in $\mathbb{R}^{n+m}$. If we assume the $n$ vectors are linear independent, up to isometry this is parametrized by the Grassmannian $G(n, m)$ of the space of $n$-dimensional subspaces in $\mathbb{R}^{n+m}$. In this paper, we investigate the geometry of Grassmannian to obtain estimates on first derivatives of the mean curvature flow. Recall the Gauss map $\gamma$ maps $(\Sigma, g_t)$ into the Grassmannian $G(n, m)$ by sending $x \in \Sigma$ to $T_x\Sigma_t \subset \mathbb{R}^{n+m}$, the tangent space of $\Sigma$ at $x$.

**Theorem A.** *The Gauss maps of a mean curvature flow $\gamma : (\Sigma, g_t) \to G(n, m)$ form a harmonic map heat flow, i.e.*

$$\frac{d}{dt}\gamma = tr\nabla d\gamma$$

*where $d\gamma$ is considered as a section of $T^*\Sigma_t \otimes \gamma^{-1}TG(n, m)$ and the trace is taken with respect to the induced metric $g_t$.*

Theorem A was conjectured by R. Hamilton and T. Ilmanen in private communications with the author. This generalizes a famous theorem of Ruh-Vilms [13] which states the Gauss map of a minimal submanifold is a harmonic map.

An immediate corollary of Theorem A is

**Corollary A.** *If $\rho$ is any convex function on $G(n,m)$, then $\gamma \circ \rho$ is a subsolution of the (nonlinear) heat equation.*

$$(\frac{d}{dt} - \Delta)\gamma \circ \rho \leq 0$$

*where $\Delta$ is the Laplace operator of the induced metric $g_t$ on $\Sigma$.*

In the stationary case, a convex function on $G(n,m)$ gives a subharmonic function on a minimal submanifold. This approach was developed by Fischer-Colbrie[7], Hildebrandt-Jost-Widman[10], Jost-Xin[11][12], Tsui-Wang[15] and [20] to obtain Bernstein type results for higher codimension minimal submanifolds .

When $\Sigma$ is compact, it follows from the maximum principle that if $\gamma(\Sigma_0)$ lies in a convex set of $G(n,m)$, so does $\gamma(\Sigma_t)$ as long as the flow exists smoothly. In the codimension-one case, the Grassmannian is a sphere and any hemisphere is a convex subset. This implies that the condition of being a graph over an affine hyperplane is preserved along the mean curvature flow, a key step in Ecker-Huisken's [4] [5] estimates for hypersurface flows.

Any $n$ form on $\mathbb{R}^{n+m}$ naturally defines a function on $G(n,m)$. We identify a class of convex functions on $G(n,m)$ in the following theorem.

**Theorem B.** *Let $\Omega$ be a simple unit $n$ form on $\mathbb{R}^{n+m}$, then $-\ln \Omega$, as a function on $G(n,m)$ is convex on a set $\Xi$ which properly contains $\{\Omega \geq \frac{1}{2}\}$.*

The set $\Xi$ will described explicitly in §3. That $-\ln \Omega$ is convex on the set $\{\Omega \geq \frac{1}{2}\}$ was proved implicitly by Jost-Xin in [11]. Instead of calculating the function $\Omega$, the authors considered the distance function on $G(n,m)$. Our approach makes the calculation simpler and more explicit, thus a stronger version is obtained.

We remark that although Theorem A and B are conceptually easier to understand, all applications rely on explicit formulae. More important information is contained in the "extra" terms. The explicit formula is derived in [19] and has been applied to study long-time existence and convergence of mean curvature flows [19] and Bernstein-type problems for minimal submanifolds in [15] and [20].

We give another example to illustrate Theorem A. On $\mathbb{R}^4$, take three orthonormal self-dual two-forms $\alpha_1, \alpha_2, \alpha_3$ and anti-self-dual two-forms $\beta_1, \beta_2, \beta_3$ as in §3 of [17]. These forms serve as coordinate functions on $G(2,2)$ under the identification
$$x \in G(2,2) \to (\alpha_i(x), \beta_i(x))$$
An element $x$ in $G(2,2)$ satisfies $\sum_{i=1}^{3}(\alpha_i(x))^2 = \sum_{i=1}^{3}(\beta_i(x))^2 = \frac{1}{2}$. Therefore $G(2,2) \simeq S_+^2(\frac{1}{\sqrt{2}}) \times S_-^2(\frac{1}{\sqrt{2}})$. It is clear that each $\alpha_i$ is a convex function on the hemisphere where $\alpha_i > 0$. Each $\alpha_i$ can be considered as a Kähler form and a surface is called symplectic with respect to $\alpha_i$ if the Gauss map lies in the region $\{\alpha_i > 0\}$. This fact is implicit in the derivation of [17] where we show being symplectic is preserved along the mean curvature flow in a Kähler-Einstein four manifold, see also Chen-Tian[3].

Indeed, the condition that the image of the Gauss map lies in the set $\{\alpha_i = 0\}$ which corresponds to a great circle on $S_+^2$ is also preserved. It is a special case of the following general theorem.

**Theorem C.** *Suppose $\Sigma_t \subset \mathbb{R}^{n+m}$, $t \in [0,T)$ are compact immersed submanifolds evolving by the mean curvature flow. If $\gamma(\Sigma_0)$ lies in a totally geodesic submanifold of $G(n,m)$, so does $\gamma(\Sigma_t)$ for $t \in [0,T)$.*

A compact Lagrangian submanifold remains Lagrangian under the mean curvature flow in a Kähler-Einstein manifold, see for example Smoczyk [14]. When the ambient space is $\mathbb{C}^n$, since the Lagrangian Grassmannian $U(n)/SO(n)$ is totally geodesic in $G(n,n)$, this also follows from Theorem C.

Theorem A and C are proved in §2. Theorem B is proved in §3. In §4, we prove $\Xi$ (see Theorem B) is a convex subset of the Grassmannian. In §5, we discuss briefly the necessary adaption when the ambient space is a general Riemannian manifold.

I am grateful to Professor S.-T. Yau for his constant support. I have benefitted from inspiring discussions with B. Andrews, R. Hamilton, T. Ilmanen, and M-P. Tsui. I would like to thank G. Huisken who encouraged me to write up this paper.

Part of this paper was completed while the author was visiting the FIM in ETH, Zürich. The author would like to thank the institute for their hospitality during his stay.

## 2 Gauss Maps

We refer to [22], [7], [10], [?] for general facts on Grassmannnian geometry. We shall adopt the description $G(n,m) = SO(n+m)/SO(n) \times SO(m)$. The tangent space of $G(n,m)$ at the identity can be identified with the space of block matrices of the form

$$\begin{bmatrix} 0 & A \\ -A^t & 0 \end{bmatrix}$$

where $A$ is an $n \times m$ matrix. The homogeneous metric on $G(n,m)$ is $ds^2 = \sum A_{i\alpha}^2$.

Let $P \in G(n,m)$ be an $n$-dimensional subspace and $T_P G = T_P G(n,m)$ be the tangent space of $G(n,m)$ at $P$. Let $e_1(s) \wedge \cdots \wedge e_n(s)$ represent a one-parameter family of $n$ planes with $\{e_1(s), \cdots, e_n(s)\}$ as their orthonormal bases so that $e_1(0) \wedge \cdots \wedge e_n(0) = e_1 \wedge \cdots \wedge e_n$ represents $P$. We have

$$\frac{d}{ds}|_{s=0} e_1(s) \wedge \cdots \wedge e_n(s) = e_1'(0) \wedge \cdots \wedge e_n + \cdots + e_1 \wedge \cdots \wedge e_n'(0).$$

By the identification $T_P G \equiv Hom(P, P^\perp)$, we may assume $e_i'(0)$ lies in the orthogonal complement $P^\perp$. The length of this tangent vector in $\wedge^n \mathbb{R}^{n+m}$ is $(\sum |e_i'(0)|^2)^{\frac{1}{2}}$. This element is identified with $\sum e_i^* \otimes e_i'(0) \in Hom(P, P^\perp)$ by the natural pairing through $e_1 \wedge \cdots \wedge e_n$. It is clearly an isometry.

Now suppose $\Sigma$ is an immersed submanifold in $\mathbb{R}^{n+m}$ with Gauss map $\gamma : \Sigma \to G(n,m)$. Thus we have a canonical identification of bundles

$$\gamma^{-1} TG \equiv T^*\Sigma \otimes N\Sigma. \tag{2.1}$$

We recall the statement of Theorem A and present the proof.

**Theorem A.** *The Gauss maps of a mean curvature flow $\gamma : (\Sigma, g_t) \to G(n,m)$ form a harmonic map heat flow, i.e.*

$$\frac{d}{dt}\gamma = tr \nabla d\gamma$$

*where $d\gamma$ is considered as a section of $T^*\Sigma_t \otimes \gamma^{-1} TG(n,m)$ and the trace is taken with respect to the induced metric $g_t$.*

*Proof.* Recall from Ruh-Vilms[13], the tension field $tr\nabla d\gamma$ of the Gauss map $\gamma$ can be identified with $\nabla H$, where $\nabla$ is the connection on the normal bundle.

We can identify $d\gamma \in \Gamma(T^*\Sigma \otimes \gamma^{-1}(TG))$ with the second fundamental form $A \in \Gamma(T^*\Sigma \otimes T^*\Sigma \otimes N\Sigma)$ through (2.1). Now $tr\nabla A = \nabla H$ by the Codazzi equation.

It suffices to show $\frac{d}{dt}\gamma = \nabla H$. Express each element in $T^*\Sigma \otimes N\Sigma$ as $\sum_i e_i^* \otimes v_i$, $v_i \in N\Sigma$, the corresponding tangent vector in the Grassmannian is

$$v_1 \wedge e_2 \wedge \cdots \wedge e_n + \cdots + e_1 \wedge \cdots \wedge e_{n-1} \wedge v_n.$$

With the identification, it is clear that the element correspond to $\nabla H$ is

$$\nabla_{e_i} H \wedge e_2 \wedge \cdots \wedge e_n + \cdots + e_1 \wedge \cdots \wedge \nabla_{e_n} H.$$

To calculate $\frac{d\gamma}{dt}$ at any space-time point $(p,t)$, we fix a coordinate system $x^1, \cdots x^n$ on $\Sigma$.

The Gauss map is determined by

$$\gamma_t = \frac{1}{\sqrt{\det g_{ij}}} \frac{\partial F_t}{\partial x^i} \wedge \cdots \wedge \frac{\partial F_t}{\partial x^n}.$$

We may assume $\frac{\partial F}{\partial x_i} = e_i$ forms an orthonormal basis at $(p,t)$. Recall the following evolution equation of the volume element:

$$\frac{d}{dt}\sqrt{\det g_{ij}} = -|H|^2 \sqrt{\det g_{ij}}.$$

On the other hand, $\frac{d}{dt}\frac{\partial F_t}{\partial x^i} = \frac{\partial H}{\partial x^i}$. By decomposition into tangent and normal parts, we have

$$\frac{\partial H}{\partial x^i} = \langle \frac{\partial H}{\partial x^i}, \frac{\partial F}{\partial x^j}\rangle g^{ik}\frac{\partial F}{\partial x^k} + \nabla_{\frac{\partial F}{\partial x^i}} H.$$

Using $\langle \frac{\partial H}{\partial x^i}, \frac{\partial F}{\partial x^j}\rangle = -\langle H, \frac{\partial^2 F}{\partial x^i \partial x^j}\rangle$, it is not hard to check

$$\frac{d}{dt}|_{t=0}\gamma_t = \nabla_{e_i} H \wedge e_2 \wedge \cdots \wedge e_n + \cdots e_1 \wedge \cdots \wedge \nabla_{e_n} H$$

and the theorem is proved.

$\square$

Now let $\rho$ be any function on $G(n,m)$, the composite function $\rho \circ \gamma : (\Sigma, g_t) \to \mathbb{R}$ satisfies the following equation

$$\frac{d}{dt}\rho \circ \gamma = d\rho(\frac{d}{dt}\gamma) = d\rho(tr\nabla d\gamma)$$

where $tr\nabla d\gamma$ is the tension field. The following calculation can be found in Proposition (2.20) of Eells-Lemaire[6]

$$\Delta_t(\rho \circ \gamma) = tr\nabla d\rho(d\gamma, d\gamma) + d\rho(tr\nabla d\gamma)$$

where $\Delta_t$ is the Laplace operator on $\Sigma_t$.

Therefore

$$(\frac{d}{dt} - \Delta)\rho \circ \gamma = -tr\nabla d\rho(d\gamma, d\gamma). \tag{2.2}$$

In case $\rho$ is convex, $\nabla d\rho$ is a positive definite quadratic form and Corollary A follows from this equation. The following theorem is a direct consequence of the maximum principle of parabolic equations.

**Theorem 2.1** *Suppose $\Sigma_t \subset \mathbb{R}^{n+m}$ are compact immersed submanifolds evolving by mean curvature flow. If $\gamma(\Sigma_0)$ lies in a convex set of $G(n,m)$, so does $\gamma(\Sigma_t)$ for $t \in [0,T)$.*

Next we recall the statement of Theorem C and present the proof.

**Theorem C.** *Suppose $\Sigma_t \subset \mathbb{R}^{n+m}$, $t \in [0,t)$ are compact immersed submanifolds evolving by the mean curvature flow. If $\gamma(\Sigma_0)$ lies in a totally geodesic submanifold of $G(n,m)$, so does $\gamma(\Sigma_t)$ for $t \in [0,T)$.*

*Proof.* Let $\mathfrak{T}$ be a totally geodesic submanifold and $d(\cdot, \mathfrak{T})$ be the distance function to $\mathfrak{T}$. We consider the second derivative of the square of the distance function, $r(\cdot) = d^2(\cdot, \mathfrak{T})$. Fix $p_0 \notin \mathfrak{T}$, we assume $p_0$ lies in a sufficiently small tubular neighborhood of $\mathfrak{T}$ without focal points. Let $\alpha_0(t), 0 \leq t \leq 1$ be a minimizing geodesic from $p_0$ to $p \in \mathfrak{T}$ which realizes the distance to $\mathfrak{T}$. We assume $\alpha_0(t)$ is parametrized so that $|\alpha'_0(t)| = d(p_0, p) = d(p_0, \mathfrak{T})$. Let $v(s), -\epsilon \leq s \leq \epsilon$ be a normal geodesic through $p_0 = v(0)$ such that $X = v'(0)$ is perpendicular to $\alpha'(0)$.

Following Proposition (3.11) of Eells-Lemaire [6] (see also Bishop-O'Neill[1]), there exists a parametrization

$$\alpha : [-\epsilon, \epsilon] \times [0, 1] \to G(n, m)$$

such that each $\alpha_s(t) = \alpha(s,t), 0 \le t \le 1$ is a minimizing geodesic from $v(s)$ to $\mathfrak{T}$ and $\alpha(s,0) = v(s)$. Denote $T = \frac{\partial \alpha}{\partial t}$ and $V = \frac{\partial \alpha}{\partial s}$, then $V$ is a Jacobi field, i.e. $\nabla_T \nabla_T V = -R(T,V)T$.

Thus $V(0,0) = v'(0)$ and $V(s,1)$ is tangent to $\mathfrak{T}$ because $\alpha(s,1)$ is contained in $\mathfrak{T}$. Also notice that

$$\langle V, T \rangle = 0 \text{ at both } p_0 = \alpha(0,0) \text{ and } p = \alpha(0,1). \tag{2.3}$$

The second variation can be calculated as in Eells-Lemaire [6] and we obtain

$$\frac{\partial^2 r}{\partial s^2}\bigg|_{s=0} = \langle \nabla_V V, T \rangle \big|_0^1 + \int_0^1 \langle \nabla_T V, \nabla_T V \rangle - \langle R(T,V)T, V \rangle \, dt. \tag{2.4}$$

Now $\langle \nabla_V V, T \rangle = 0$ at $p_0 = \alpha(0,0)$ because $v(s)$ is a geodesic. At $p = \alpha(0,1)$, $T$ is normal to $\mathfrak{T}$ and $V$ is tangent to $\mathfrak{T}$, therefore $\langle \nabla_V V, T \rangle$ represents a second fundamental form of $\mathfrak{T}$; this terms vanishes as well because $\mathfrak{T}$ is totally geodesic.

We claim $\langle V, T \rangle = 0$ along $\alpha_0$. Indeed, taking derivative of $\langle V, T \rangle$ with respect to $t$ twice, we obtain

$$T \langle \nabla_T V, T \rangle = \langle \nabla_T \nabla_T V, T \rangle = -\langle R(T,V)T, T \rangle = 0.$$

Thus $\langle \nabla_T V, T \rangle = T \langle V, T \rangle$ is a constant function in $t$, or $\langle V, T \rangle$ is a linear function in $t$ along $\alpha_0$. Since $\langle V, T \rangle = 0$ at $p_0$ and $p$, we conclude $\langle V, T \rangle = 0$ along $\alpha_0$.

From $\langle V, T \rangle = 0$, we deduce $\langle R(T,V)T, V \rangle \le K_1 |V|^2 |T|^2$ along $\alpha_0$, where $K_1$ is an upper bound of the sectional curvature of $G(n,m)$. Moreover, along $\alpha_0$ we have $|V|^2 \le K_2 |X|^2$ for a constant $K_2$ that depends on the sectional curvature and the size of the tubular neighborhood by a comparison argument. Now (2.4) implies

$$\frac{\partial^2 r}{\partial s^2}\bigg|_{s=0} \ge -K_3 |X|^2$$

or

$$\nabla dr(X,X) \ge -K_3 |X|^2$$

at any point $p_0$ in a sufficiently small tubular neighborhood of $\mathfrak{T}$.

Combine this equation with equation (2.2), we derive

$$(\frac{d}{dt} - \Delta)(r \circ \gamma) \leq K' r \circ \gamma. \tag{2.5}$$

$K'$ involves the second fundamental form (recall $A = d\gamma$), but as long as the flow exists smoothly $K'$ is bounded. The assumption implies $r \circ \gamma = 0$ initially and this remains true afterwards by applying the maximum principle to (2.5). □

## 3 Applications

Note the summation convention, repeated indices are summed over, is adopted in the rest of this article. We first define the set $\Xi$ in the statement of Theorem B. $\Omega$, as a simple unit $n$-form, is dual to an $n$-subspace $Q$. Given any $P$ that can be written as the graph of a linear transformation $L_P : Q \to Q^\perp$ over $Q$. Denote by $\lambda_i(P)$ the singular values $L_P$ or the eigenvalues of the symmetric matrix $\sqrt{L_P^T L_P}$. $\Xi$ is then defined by

$$\Xi = \{P \in G(n,m) \mid P \text{ is a graph over } Q \text{ and } |\lambda_i \lambda_j| \leq 1 \text{ for any } i \neq j.\} \tag{3.1}$$

It is not hard to see $\Omega(P) = \frac{1}{\sqrt{\prod(1+\lambda_i^2)}}$. Therefore

$$\Omega \geq \frac{1}{2} \text{ implies } |\lambda_i \lambda_j| \leq 1.$$

**Theorem B.** *Let $\Omega$ be a simple unit $n$ form on $\mathbb{R}^{n+m}$, then $-\ln \Omega$, as a function on $G(n,m)$ is convex on $\Xi$ which properly contains the set $\{\Omega \geq \frac{1}{2}\}$.*

*Proof.* Given $P = e_1 \wedge \cdots \wedge e_n \in G(n,m)$ and suppose the orthogonal complement $P^\perp$ is spanned by $\{e_{n+\alpha}\}_{\alpha=1\cdots m}$. By a formula of Wong[22], a geodesic through $P$ parametrized by arc length is given as $P_s$ spanned by $\{e_i + z_{i\alpha}(s) e_{n+\alpha}\}_{i=1\cdots n}$ such that $Z = [z_{i\alpha}(s)]$ is a $n \times m$ matrix which satisfies the following ordinary differential equation:

$$Z'' - 2Z' Z^T (I + ZZ^T)^{-1} Z' = 0.$$

We assume $z_{i\alpha}(0) = 0$, i.e. $P_0 = P$, then the geodesic equation implies $z''_{i\alpha}(0) = 0$. We also denote $z'_{i\alpha}(0) = \mu_{i\alpha}$ in the following calculation.

We shall calculate the second derivative of the ln of

$$\mathfrak{p}(s) = \Omega(P_s) = \frac{1}{\sqrt{\det g_{ij}}} \Omega(e_1 + z_{1\alpha}(s)e_{n+\alpha}, \cdots, e_n + z_{n\alpha}(s)e_{n+\alpha}) \quad (3.2)$$

where $g_{ij}(s) = \delta_{ij} + z_{i\alpha}(s)z_{j\alpha}(s)$. By direct calculation,

$$g'_{ij}(0) = 0, \text{ and } g''_{ij}(0) = 2\mu_{i\alpha}\mu_{j\alpha}.$$

Therefore $g = \det g_{ij}$ satisfies

$$(\frac{1}{\sqrt{g}})'(0) = 0, \text{ and } (\frac{1}{\sqrt{g}})''(0) = -\mu_{i\alpha}^2. \quad (3.3)$$

Differentiate equation (3.2) and plug in (3.3), we obtain

$$\mathfrak{p}'(0) = \mu_{1\alpha}\Omega(e_{n+\alpha}, e_2, \cdots, e_n) + \cdots + \mu_{n\alpha}\Omega(e_1, \cdots, e_{n-1}, e_{n+\alpha})$$
$$\mathfrak{p}''(0) = -(\sum_{i,\alpha} \mu_{i\alpha}^2)\Omega(e_1, e_2, \cdots, e_n) + 2[\mu_{1\alpha}\mu_{2\beta}\Omega(e_{n+\alpha}, e_{n+\beta}, \cdots, e_n) + \cdots],$$

the expression in the bracket runs through $(i,j)$ with $i < j$, the general form is $\mu_{i\alpha}\mu_{j\beta}\Omega(e_1, \cdots, e_{n+\alpha}, \cdots, e_{n+\beta}, \cdots, e_n)$ where we replace $e_i$ and $e_j$ by $e_{n+\alpha}$ and $e_{n+\beta}$, respectively. Now

$$(\ln \mathfrak{p})'' = \frac{1}{\mathfrak{p}^2}[\mathfrak{p}''\mathfrak{p} - (\mathfrak{p}')^2]$$

We shall assume $\Omega(P) > 0$, therefore $P$ can be written as a graph over the plane dual to $\Omega$. By singular value decomposition, we can choose the basis $e_1, \cdots e_n$ for $P$ and $e_{n+1}, \cdots e_{n+m}$ for $P^\perp$ such that

$$e_i = \frac{1}{\sqrt{1+\lambda_i^2}}(a_i + \lambda_i a_{n+i}) \text{ and } e_{n+\alpha} = \frac{1}{\sqrt{1+\lambda_\alpha^2}}(a_{n+\alpha} - \lambda_\alpha a_\alpha) \quad (3.4)$$

for $i = 1 \cdots n$ and $\alpha = 1 \cdots m$ where we pretend $\lambda_i = 0$ and $\lambda_\alpha = 0$ for $i, \alpha > \min\{m, n\}$. Here $\{a_i\}$ is an orthonormal basis for the plane dual to $\Omega$ and $\{a_{n+\alpha}\}$ an orthonormal basis for the orthogonal complement.

$$\mathfrak{p}'(0) = -\mu_{i,n+i}\lambda_i \mathfrak{p}(0)$$

$$\mathfrak{p}''(0) = [-(\sum_{i,\alpha}\mu_{i,n+\alpha}^2) + 2\sum_{i<j}\mu_{i,n+i}\mu_{j,n+j}\lambda_i\lambda_j - 2\sum_{i<j}\mu_{i,n+j}\mu_{j,n+i}\lambda_i\lambda_j]\mathfrak{p}(0)$$

Therefore

$$(\ln \mathfrak{p})''(0) = -(\sum_{i,n+\alpha}\mu_{i,n+\alpha}^2) + 2\sum_{i<j}\mu_{i,n+i}\mu_{j,n+j}\lambda_i\lambda_j - 2\sum_{i<j}\mu_{i,n+j}\mu_{j,n+i}\lambda_i\lambda_j - (\sum_i \mu_{i,n+i}\lambda_i)^2$$
$$= -(\sum_{i,\alpha}\mu_{i,n+\alpha}^2) - 2\sum_{i<j}\mu_{i,n+j}\mu_{j,n+i}\lambda_i\lambda_j - \sum_i(\mu_{i,n+i}\lambda_i)^2.$$
(3.5)

By completing square we derive $-\ln \mathfrak{p}$ is a convex function of $s$ if $|\lambda_i\lambda_j| \leq 1$ for any $i \neq j$. Since we can perform this calculation in any direction, $-\ln \Omega$ is a convex function on $\Xi$. □

Theorem B should be compared with the explicit formula derived in [19].

## 4 Grassmannian Convexity of $\Xi$

First we give a new characterization of $\Xi$ in terms of the positivity of a bilinear form. Consider the bilinear form $S$ on $\mathbb{R}^{n+m}$ by

$$S(X,Y) = \langle \pi_1(X), \pi_1(Y) \rangle - \langle \pi_2(X), \pi_2(Y) \rangle \text{ for } X, Y \in \mathbb{R}^{n+m}$$

where $\pi_1$ and $\pi_2$ are projections from $\mathbb{R}^{n+m} = \mathbb{R}^n \oplus \mathbb{R}^m$ to $\mathbb{R}^n$ and to $\mathbb{R}^m$, respectively. Let $P$ be an $n$-subspace of $\mathbb{R}^{n+m}$ that is a graph over the base $Q \cong \mathbb{R}^n$; we denote the restriction of $S$ to $P$ by $S|_P$.

For any (not necessarily orthonormal) basis $\{e_i\}$ of $P$, $S|_P$ is represented by

$$S_{ij} = S(e_i, e_j).$$

We also define

$$g_{ij} = \langle e_i, e_j \rangle \text{ and } \sigma_i^j = \frac{1}{2}(g^{ik}S_{kj} + g^{jk}S_{ki})$$

where $g^{ij}$ is the inverse to $g_{ij}$. Then $\sigma = (\sigma_i^j) : P \to P$ becomes a self-adjoint map and satisfies
$$\langle \sigma(P)(X), Y \rangle = S|_P(X, Y).$$

$\sigma$ induces a linear map on $\wedge^2 P$ by
$$\sigma(X \wedge Y) = \sigma(X) \wedge Y + X \wedge \sigma(Y).$$

The collection $\{e_i \wedge e_j\}_{i<j}$ forms a basis for $\wedge^2 P$ and in terms of this basis

$$\sigma(e_i \wedge e_j) = \sum_{k,l} \sigma_i^k e_k \wedge e_j + \sigma_j^l e_i \wedge e_l = \sum_{k<l}(\sigma_i^k \delta_j^l + \sigma_j^l \delta_i^k - \sigma_i^l \delta_j^k - \sigma_j^k \delta_i^l) e_k \wedge e_l.$$

We can use $\sigma$ to characterize $\Xi$, in fact
$$\Xi = \{P \in G(n,m) \mid \min_{\omega \in \wedge^2 P} \langle \sigma(P)(\omega), \omega \rangle \geq 0\}.$$

To see this, apply singular value decomposition to find an orthonormal basis as in the last section $\{e_i = \frac{1}{\sqrt{1+\lambda_i^2}}(a_i + \lambda_i a_{n+i})\}$ for $P$, then
$$S_{ij} = \sigma_j^i = \frac{1 - \lambda_i^2}{1 + \lambda_i^2} \delta_{ij}.$$

The coefficient of $\sigma$ on $\wedge^2 P$ is

$$\sigma_{(ij)}^{(kl)} = \sigma_i^k \delta_j^l + \sigma_j^l \delta_i^k = \frac{1 - \lambda_i^2}{1 + \lambda_i^2} \delta_i^k \delta_j^l + \frac{1 - \lambda_j^2}{1 + \lambda_j^2} \delta_j^l \delta_i^k = \frac{2(1 - \lambda_i^2 \lambda_j^2)}{(1 + \lambda_i^2)(1 + \lambda_j^2)} \delta_i^k \delta_j^l \quad (4.1)$$

if $i < j$ and $k < l$. Therefore $\sigma(P)$ being positive definite on $\wedge^2 P$ is the same as the area-decreasing condition $|\lambda_i \lambda_j| \leq 1$.

Now we use this characterization by $\sigma$ to prove the convexity of $\Xi$.

**Theorem 4.1** $\Xi$ *is a convex subset with respect to the Grassmannian metric.*

*Proof.* The proof is inspired by Hamilton's maximum principle [8] [9] for tensors. Let $P$ be a boundary point of $\Xi$, so $\sigma(P)$ is non-negative definite on $\wedge^2 P$. Let $\omega \in \wedge^2 P$ be a zero eigenvector of $\sigma(P)$ so that $\langle \sigma(P)\omega, \omega \rangle = 0$. Consider a (Grassmannian) geodesic $P(s)$ through $P$ and an extension of $\omega$, $\omega_s$ on $P(s)$ and denote $f(s) = \langle \sigma(P(s))(\omega_s), \omega_s \rangle$. To check the convexity, it suffices to show for any geodesic $P(s)$ we can find an (arbitrary) extension $\omega_s$ so that $f'(0) = 0$ and $f''(0) < 0$. We remark that an arbitrary extension of $\omega$ is good enough as the minimum function $\min_{\omega \in \wedge^2 P}\langle \sigma(P)(\omega), \omega\rangle$ is always less than or equal to $f(s)$ along $P(s)$.

As in the previous section, we choose an orthonormal basis $\{e_i\}$ for $P$, so that $P = e_1 \wedge \cdots \wedge e_n \in G(n,m)$ and suppose the orthogonal complement $P^\perp$ is spanned by $\{e_{n+\alpha}\}_{\alpha=1\cdots m}$. A geodesic parametrized by arc length is given as $P_s$ spanned by $\{e_i + z_{i\alpha}(s)e_{n+\alpha}\}_{i=1\cdots n}$.

Denote
$$S_{ij}(s) = S(e_i + z_{i\alpha}(s)e_{n+\alpha}, e_j + z_{j\beta}(s)e_{n+\beta}), \qquad (4.2)$$
$$g_{ij}(s) = \langle e_i + z_{i\alpha}(s)e_{n+\alpha}, e_j + z_{j\beta}(s)e_{n+\beta}\rangle = \delta_{ij} + z_{i\alpha}(s)z_{j\alpha}(s).$$

For any element of $\wedge^2 P_s$,
$$\omega_s = \sum_{i<j} \omega^{ij}(s)(e_i + z_{i\alpha}(s)e_{n+\alpha}) \wedge (e_j + z_{j\beta}(s)e_{n+\beta}),$$

we have
$$|\omega_s|^2 = \sum_{i<j,k<l} \omega^{ij}(s)\omega^{kl}(s)(g_{ik}(s)g_{jl}(s) - g_{il}(s)g_{jk}(s)). \qquad (4.3)$$

We shall choose $\omega_s$ so that $|\omega_s|^2$ is constant up to second order at $s = 0$. On the other hand,
$$f(s) = \langle \sigma(P_s)(\omega_s), \omega_s\rangle = \sum_{i<j,k<l} \omega^{ij}(s)\omega^{kl}(s)\sigma_{(ij)(kl)}(s)$$

where
$$\sigma_{(ij)(kl)}(s) = S_{ik}(s)g_{jl}(s) + S_{jl}(s)g_{ik}(s) - S_{il}(s)g_{jk}(s) - S_{jk}(s)g_{il}(s). \qquad (4.4)$$

We recall that $z_{i\alpha}(0) = z''_{i\alpha}(0) = 0$ and $z'_{i\alpha}(0) = \mu_{i\alpha}$, also
$$g'_{ij}(0) = 0, \text{ and } g''_{ij}(0) = 2\mu_{i\alpha}\mu_{j\alpha}.$$

In the following calculations, all derivatives are taken at $s = 0$. By equation (4.2),
$$S'_{ij} = \mu_{i\alpha}S(e_{n+\alpha}, e_j) + \mu_{j\beta}S(e_i, e_{n+\beta}), \quad S''_{ij} = 2\mu_{i\alpha}\mu_{j\beta}S(e_{n+\alpha}, e_{n+\beta}),$$

By (4.4), we derive for $i < j$ and $k < l$,

$$\sigma'_{(ij)(kl)} = S'_{ik}\delta_{jl} + S'_{jl}\delta_{ik} - S'_{il}\delta_{jk} - S'_{jk}\delta_{il}$$

and

$$(\sigma_{(ij)(kl)})'' = S''_{ik}\delta_{jl} + S''_{jl}\delta_{ik} - S''_{il}\delta_{jk} - S''_{jk}\delta_{il} + S_{ik}g''_{jl} + S_{jl}g''_{ik} - S_{il}g''_{jk} - S_{jk}g''_{il}.$$

Since $\omega$ is a zero eigenvector of $\sigma$ on $\wedge^2 P$, $\omega = e_i \wedge e_j$ where $e_i$ and $e_j$ are eigenvectors of $\sigma$ on $P$. By reordering the basis, we may assume $\omega = e_1 \wedge e_2$. We extend $\omega$ to $\omega_s$ so that $(\omega^{ij})'(0) = 0$ to make $(|\omega_s|^2)'(0) = 0$.

Now

$$f'(0) = (\sigma_{(ij)(kl)})'\omega^{ij}\omega^{kl} = (\sigma_{(12)(12)})' = 2\mu_{1\alpha}S(e_{n+\alpha}, e_1) + 2\mu_{2\alpha}S(e_{n+\alpha}, e_2)$$

By (3.4), $S(e_{n+\alpha}, e_k) = -\frac{2\lambda_k}{1+\lambda_k^2}\delta_{\alpha k}$, so $f'(0) = 0$ implies

$$\frac{\lambda_1}{1+\lambda_1^2}\mu_{11} + \frac{\lambda_2}{1+\lambda_2^2}\mu_{22} = 0. \tag{4.5}$$

On the other hand, $f(0) = 0$ implies $\sigma^{(12)}_{(12)} = 0$ or $\lambda_1\lambda_2 = 1$. It follows that $\frac{\lambda_1}{1+\lambda_1^2} = \frac{\lambda_2}{1+\lambda_2^2}$ and by (4.5),

$$\mu_{11}^2 = \mu_{22}^2. \tag{4.6}$$

To keep $(|\omega_s|^2)''(0) = 0$, it suffices to set

$$2(\omega^{12})''(0) = -g''_{11} - g''_{22} = -2\sum_\alpha \mu_{1\alpha}^2 - 2\sum_\alpha \mu_{2\alpha^2}.$$

We shall assume $m \geq n$ in the following calculation, the case $m < n$ can be carried out similarly. The second derivative of $f$ can be calculated in the following.

$$\begin{aligned}-(f)''(0) &= -(\sigma_{(ij)(kl)})''\omega_{ij}\omega_{kl} - (\sigma_{(ij)(kl)})(\omega_{ij})''\omega_{kl} - (\sigma_{(ij)(kl)})\omega_{ij}(\omega_{kl})'' \\ &= -(\sigma_{(12)(12)})'' - 2(\omega^{12})''\sigma_{(12)(12)} \\ &= -2\mu_{1\alpha}\mu_{1\beta}S(e_{n+\alpha}, e_{n+\beta}) - 2\mu_{2\alpha}\mu_{2\beta}S(e_{n+\alpha}, e_{n+\beta}) - 2(\sum_\alpha \mu_{1\alpha}^2)S_{22} - 2(\sum_\alpha \mu_{2\alpha}^2)S_{11} \\ &\quad + 2(\sum_\alpha \mu_{1\alpha}^2 + \mu_{2\alpha}^2)(S_{11} + S_{22})\end{aligned}$$

Write $T_{\alpha\alpha} = S(e_{n+\alpha}, e_{n+\alpha})$ and deduce from (3.4) $T_{11} = -S_{11}$ and $T_{22} = -S_{22}$, the last expression is equal to

$$2\sum_\alpha \mu_{1\alpha}^2 S_{11} + 2\sum_\alpha \mu_{2\alpha}^2 S_{22} + 2\sum_k \mu_{1\alpha}^2 T_{\alpha\alpha} + 2\sum_k \mu_{2\alpha}^2 T_{\alpha\alpha}$$
$$= 4\mu_{11}^2 S_{11} + 4\mu_{22}^2 S_{22} + 2\mu_{12}^2(S_{11} + S_{22}) + 2\mu_{21}^2(S_{11} + S_{22})$$
$$+ 2\sum_{\alpha \geq 3} \mu_{1\alpha}^2(S_{11} - T_{\alpha\alpha}) + 2\sum_{\alpha \geq 3} \mu_{2\alpha}^2(S_{22} - T_{\alpha\alpha}).$$

Now $S_{11} + S_{22} = 0$ and it is not hard to check $S_{11} - T_{\alpha\alpha}$ and $S_{22} - T_{\alpha\alpha}$ are both non-negative for any $\alpha \geq 3$. Therefore

$$-f''(0) \geq 4\mu_{11}^2 S_{11} + 4\mu_{22}^2 S_{22}.$$

It follows from equation (4.6) that $f''(0)$ is non-positive and the theorem is proved. $\square$

We remark a similar calculation in the mean curvature flow case is presented in Tsui-Wang [16].

## 5 Riemannian Ambient Manifolds

When applying theorem A and C to general ambient Riemannian submanifold $M$, the validity of the theorem depends on the curvature of $M$. The Gauss map $\gamma$ is a section of the Grassmannian bundle $\mathfrak{G}$ over the ambient space $M$. Now $\gamma^{-1}(T\mathfrak{G})$ can be identified with $T^*\Sigma \otimes N\Sigma$ as before. $d\gamma$ is a section of $T^*\Sigma \otimes T^*\Sigma \otimes N\Sigma$ which is exactly the second fundamental form $A$.

By the Codazzi equation,

$$tr\nabla d\gamma = tr\nabla A = \nabla H + (R(\cdot, e_k)e_k)^\perp \tag{5.1}$$

where $R(X, Y)Z = -\nabla_X \nabla_Y Z + \nabla_Y \nabla_X Z - \nabla_{[X,Y]} Z$ is the curvature operator of $M$ and $\{e_k\}_{k=1\cdots n}$ is an orthonormal basis for $T\Sigma$. Notice that our convention is $\langle R(X,Y)X, Y \rangle > 0$ if $M$ has positive sectional curvature. (5.1) is considered as an equation of sections of $T^*\Sigma \otimes N\Sigma$.

Therefore, as in Theorem A we have

$$\frac{d}{dt}\gamma = tr\nabla d\gamma + (R(e_k, \cdot)e_k)^\perp$$

Now if $\Omega$ is any $n$ form on $M$, it defines a function on the Grassmannian bundle. When $\Omega = \rho$ is a parallel form, we get as before

$$(\frac{d}{dt} - \Delta)\rho \circ \gamma = -tr\nabla d\rho(d\gamma, d\gamma) + \nabla\rho \cdot (R(e_k, \cdot)e_k)^\perp \qquad (5.2)$$

where $\nabla\rho$ is consider as a tangent vector on the Grassmannian. The curvature term usually prefers positive ambient curvature, we refer to [17], [18] and [19] for explicit calculations under various curvature conditions.

When $\Omega$ is not parallel, the formula involves the covariant derivatives of $\Omega$ and the general equation was derived in [21].